\documentclass[twoside]{article}
\usepackage[accepted]{aistats2012}

\usepackage{epsf}
\usepackage{graphicx}
\usepackage{amsmath}
\usepackage{amssymb}
\usepackage{algorithmic}
\usepackage{algorithm}
\usepackage{epsfig}
\usepackage{natbib}
\bibpunct{(}{)}{;}{a}{,}{,}

\newcommand{\B}{\mathcal B}
\newcommand{\Er}{\mathcal{E}_{xploring}}

\newcommand{\A}{\mathcal A}
\newcommand{\Alg}{\mathcal Alg}

\newcommand{\D}{\mathcal B}

\renewcommand{\P}{\mathbb P}
\newcommand{\R}{\mathbb R}

\newlength{\minipagewidth}
\newcommand{\bookbox}[1]{\small
\par\medskip\noindent
\framebox[\columnwidth]{
\begin{minipage}{\columnwidth-0.5cm} {#1} \end{minipage} } \par\medskip }

\newcommand{\beq}{\begin{equation}}
\newcommand{\eeq}{\end{equation}}

\newcommand{\beqa}{\begin{eqnarray}}
\newcommand{\eeqa}{\end{eqnarray}}

\newcommand{\beqan}{\begin{eqnarray*}}
\newcommand{\eeqan}{\end{eqnarray*}}

\newcommand{\E}{\mathbb{E}}

\newcommand{\si}{\sigma}

\newcommand{\ind}[1]{\mathbb I\left\lbrace {#1} \right\rbrace}

\newtheorem{assumption}{Assumption}
\newtheorem{theorem}{Theorem}

\begin{document}

%

%

\runningtitle{Bandit Theory meets Compressed Sensing for high-dimensional Stochastic Linear Bandit}
\twocolumn[

\aistatstitle{Bandit Theory meets Compressed Sensing for high-dimensional Stochastic Linear Bandit}

\aistatsauthor{ Alexandra Carpentier \And Remi Munos }

\aistatsaddress{ Sequel team, INRIA Lille - Nord Europe \And Sequel team, INRIA Lille - Nord Europe} ]

\begin{abstract}
  We consider a linear stochastic bandit problem where the dimension $K$ of the unknown parameter $\theta$ is larger than the sampling budget $n$. Since usual linear bandit algorithms have a regret of order $O(K\sqrt{n})$, it is in general impossible to obtain a sub-linear regret without further assumption. In this paper we make the assumption that $\theta$ is $S-$sparse, i.e.~has at most $S-$non-zero components, and that the set of arms is the unit ball for the $||.||_2$ norm. We combine ideas from Compressed Sensing and Bandit Theory to derive an algorithm with a regret bound in $O(S\sqrt{n})$.
We detail an application to the problem of optimizing a function that depends on many variables but among which only a small number of them (initially unknown) are relevant.
\end{abstract}

\vspace{-0.2cm}

\section*{Introduction}

We consider a linear stochastic bandit problem in high dimension $K$. At each round $t$, from $1$ to $n$, the player chooses an arm $x_t$ in a fixed set of arms and receives a reward $r_t= \langle x_t, \theta + \eta_t \rangle $, where $\theta \in \R^K$ is an unknown parameter and $\eta_t$ is a noise term. Note that $r_t$ is a (noisy) projection of $\theta$ on $x_t$. The goal of the learner is to maximize the sum of rewards.

We are interested in cases where the number of rounds is much smaller than the dimension of the parameter, i.e.~$n \ll K$. This is new in bandit literature but useful in practice, as illustrated by the problem of gradient ascent for a high-dimensional function, described later.

In this setting it is in general impossible to estimate $\theta$ in an accurate way (since there is not even one sample per dimension). It is thus necessary to restrict the setting, and the assumption we consider here is that $\theta$ is $S$-\textit{sparse} (i.e., at most $S$ components of $\theta$ are non-zero). We assume also that the set of arms to which $x_t$ belongs is the unit ball with respect to the $||.||_2$ norm, induced by the inner product.


%

\vspace{-0.2cm}

\paragraph{Bandit Theory meets Compressed Sensing} 

This problem poses the fundamental question at the heart of bandit theory, namely the exploration\footnote{Exploring all directions enables to build a good estimate of all the components of $\theta$ in order to deduce which arms are the best.} versus exploitation\footnote{Pulling the empirical best arms in order to maximize the sum of rewards.} dilemma. Usually, when the dimension $K$ of the space is smaller than the budget $n$, it is possible to project the parameter $\theta$ \textit{at least once} on \textit{each} directions of a basis (e.g.~the canonical basis) which enables to explore efficiently. However, in our setting where $K\gg n$, this is not possible anymore, and we use the sparsity assumption on $\theta$ to build a clever exploration strategy.


\textit{Compressed Sensing} (see e.g.~\citep{candes2007dantzig,chen1999atomic,blumensath2009iterative}) provides us with a exploration technique that enables to estimate $\theta$, or more simply its support, \textit{provided that $\theta$ is sparse}, with few measurements. The idea is to project $\theta$ on random (isotropic) directions $x_t$ such that each reward sample provides equal information about \textit{all} coordinates of $\theta$. This is the reason why we choose the set of arm to be the unit ball. Then, using a regularization method (Hard Thresholding, Lasso, Dantzig selector...), one can recover the support of the parameter. 
Note that although Compressed Sensing enables to build a good estimate of $\theta$, it is not designed for the purpose of maximizing the sum of rewards. Indeed, this exploration strategy is uniform and non-adaptive (i.e., the sampling direction $x_t$ at time $t$ does not depend on the previously observed rewards $r_1,\dots,r_{t-1}$). 

On the contrary, \textit{Linear Bandit Theory} (see e.g.~\citet{rusmevichientong2008linearly,dani2008stochastic,filippiparametric} and the recent work by \citet{Csaba}) addresses this issue of maximizing the sum of rewards by efficiently balancing between exploration and exploitation. The main idea of our algorithm is to use Compressed Sensing to estimate the (small) support of $\theta$, and combine this with a linear bandit algorithm with a set of arms restricted to the estimated support of $\theta$. 




{\bf Our contributions} are the following:

\vspace{-0.2cm}
\begin{itemize}
 \item We provide an algorithm, called SL-UCB (for Sparse Linear Upper Confidence Bound) that mixes ideas of Compressed Sensing and Bandit Theory and provide a regret bound\footnote{We define the notion of regret in Section \ref{s:setting}.} of order $O(S\sqrt{n})$.
\item We detailed an application of this setting to the problem of gradient ascent of a high-dimensional function that depends on a small number of relevant variables only (i.e., its gradient is sparse). We explain why the setting of gradient ascent can be seen as a bandit problem and report numerical experiments showing the efficiency of SL-UCB for this high-dimensional optimization problem.
\end{itemize}

\vspace{-0.2cm}

The topic of sparse linear bandits is also considered in the paper~\citep{Yasin} published simultaneously. Their regret bound scales as $O(\sqrt{K S n})$ (whereas ours do not show any dependence on $K$) but they do not make the assumption that the set of arms is the Euclidean ball and their noise model is different from ours.

In Section \ref{s:setting} we describe our setting and recall a result on linear bandits. Then in Section~\ref{s:SL-UCB} we describe the SL-UCB algorithm and provide the main result. In Section~\ref{s:grad.asc} we detail the application to gradient ascent and provide numerical experiments. 

\vspace{-0.2cm}

\section{Setting and a useful existing result} \label{s:setting}
\vspace{-0.1cm}

\subsection{Description of the problem}

We consider a linear bandit problem in dimension $K$. An algorithm (or strategy) $\Alg$ is given a budget of $n$ pulls. At each round $1\leq t \leq n$ it selects an arm $x_t$ in the set of arms $\D_{K}$, which is the unit ball for the $||.||_2$-norm induced by the inner product. It then receives a reward
\begin{equation*}
r_t = \langle x_t,\theta+\eta_t\rangle ,
\end{equation*}
where $\eta_t\in \R^K$ is an i.i.d.~white noise\footnote{This means that $\E_{\eta_t}(\eta_{k,t}) = 0$ for every $(k,t)$, that the $(\eta_{k,t})_k$ are independent and that the $(\eta_{k,t})_t$ are i.i.d..} that is independent from the past actions, i.e.~from $\Big\{(x_{t'})_{t'\leq t}\Big\}$, and $\theta \in \R^K$ is an unknown parameter.

We define the \textit{performance} of algorithm $\Alg$ as

\vspace{-0.2cm}
\begin{equation}\label{def:performance}
L_n(\Alg) = \sum_{t=1}^n  \langle \theta,x_t\rangle .
\end{equation}

Note that $L_n(\Alg)$ differs from the sum of rewards $\sum_{t=1}^n r_t$ but is close (up to a $O(\sqrt{n})$ term) in high probability. Indeed, $\sum_{t=1}^n \langle \eta_t,x_t\rangle$ is a Martingale, thus if we assume that the noise $\eta_{k,t}$ is bounded by $\frac{1}{2}\si_k$ (note that this can be extended to sub-Gaussian noise), Azuma's inequality implies that with probability $1-\delta$, we have $\sum_{t=1}^n r_t = L_n(\Alg)+ \sum_{t=1}^n \langle \eta_t,x_t\rangle \leq L_n(\Alg) + \sqrt{2\log(1/\delta)} ||\si||_2 \sqrt{n}$.


If the parameter $\theta$ were known, the best strategy $\Alg^*$ would always pick $x^* = \arg\max_{x \in \D_{K}} \langle \theta,x\rangle = \frac{\theta}{||\theta||_2}$ and obtain the performance:

\vspace{-0.3cm}
\begin{equation}\label{def:performance*}
L_n(\Alg^*) = n ||\theta||_2.
\end{equation}

We define the \textit{regret} of an algorithm $\Alg$ with respect to this optimal strategy as

\vspace{-0.3cm}
\begin{equation}\label{def:regret}
R_n(\Alg) = L_n(\Alg^*) - L_n(\Alg).
\end{equation}


We consider the class of algorithms that do not know the parameter $\theta$. Our objective is to find an adaptive strategy $\Alg$ (i.e.~that makes use of the history $\{(x_1,r_1),\ldots,(x_{t-1},r_{t-1})\}$ at time $t$ to choose the next state $x_t$) with smallest possible regret.

For a given $t$, we write $X_t = (x_1;\ldots;x_t)$ the matrix in $\R^{K\times t}$ of all chosen arms, and $R_t = (r_1, \ldots, r_t)^T$ the vector in $\R^t$ of all rewards, up to time $t$.

In this paper, we consider the case where the dimension $K$ is much larger than the budget, i.e., $n \ll K$.
As already mentioned, in general it is impossible to estimate accurately the parameter and thus achieve a sub-linear regret. This is the reason why we make the assumption that $\theta$ is $S-$sparse with $S<n$.

\subsection{A useful algorithm for Linear Bandits}\label{ss:dani}

\begin{figure*}[!hbtp]
\begin{minipage}{0.5\textwidth}
\bookbox{
\begin{algorithmic}
\STATE \hspace{5mm}\textbf{Input:} $\D_{d}$, $\delta$
\STATE \hspace{5mm}\textbf{Initialization:}
\STATE \hspace{5mm}$A_1 = I_d$, $\hat{\theta}_1 = 0$, $\beta_t = 128 d  (\log(n^2/\delta))^2$.
\STATE \hspace{5mm}\textbf{for} $t = 1,\ldots, n$ \textbf{do}
\STATE \hspace{10mm}Define $B_t = \{\nu: ||\nu - \hat{\theta_t}||_{2,A_t} \leq  \sqrt{\beta_t} \}$
\STATE \hspace{10mm}Play $x_t = \arg\max_{x \in \D_{d}} \max_{\nu \in B_t} \langle \nu,x\rangle $.
\STATE \hspace{10mm}Observe $r_t = \langle x_t,\theta+\eta_t\rangle $.
\STATE \hspace{10mm}Set $A_{t+1} = A_t + x_t x_t'$, $\hat{\theta}_{t+1} = A_{t+1}^{-1} X_t R_t$.
 \STATE \hspace{5mm}\textbf{end for}
\end{algorithmic}}
\end{minipage}
\hspace{0.5cm}
\begin{minipage}{0.5\textwidth}
\includegraphics[width=7cm]{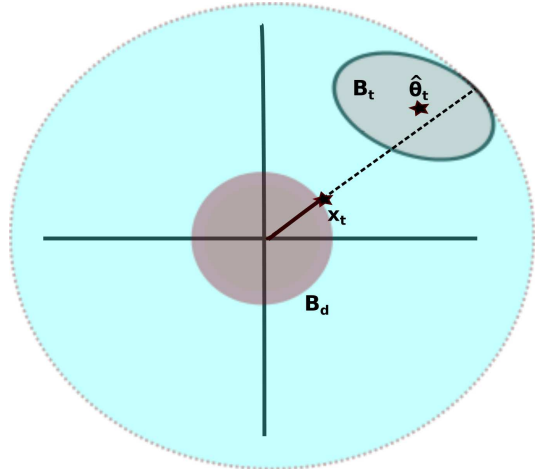}
\end{minipage}
\caption{Algorithm $ConfidenceBall_2$ ($CB_2$) adapted for an action set of the form $\D_{d}$ (Left), and illustration of the maximization problem that defines $x_t$ (Right).}
\label{fig:alglinban1}
\end{figure*}

We now recall the algorithm $ConfidenceBall_2$ (abbreviate by $CB_2$) introduced in \citet{dani2008stochastic} and mention the corresponding regret bound. $CB_2$ will be later used in the SL-UCB algorithm described in the next Section to the subspace restricted to the estimated support of the parameter.

This algorithm is designed for stochastic linear bandit in dimension $d$ (i.e.~the parameter $\theta$ is in $\R^d$) where $d$ is \textit{smaller} than the budget $n$. 



The pseudo-code of the algorithm is presented in Figure~\ref{fig:alglinban1}. The idea is to build an ellipsoid of confidence for the parameter $\theta$, namely $B_t = \{\nu: ||\nu - \hat{\theta_t}||_{2,A_t} \leq  \sqrt{\beta_t} \}$ where $||u||_{2,A}=u^TAu$ and $\hat{\theta}_{t} = A_{t}^{-1} X_{t-1} R_{t-1}$, and to pull the arm with largest inner product with a vector in $B_t$, i.e.~the arm $x_t = \arg\max_{x \in \D_{d}} \max_{\nu \in B_t} \langle \nu,x\rangle $.

Note that this algorithm is intended for general shapes of the set of arms. We can thus apply it in the particular case where the set of arms is the unit ball $\D_{d}$ for the $||.||_2$ norm in $\R^d$. This specific set of arms is simpler for two reasons. First, it is easy to define a span of the set of arms since we can simply choose the canonical basis of $\R^d$. Then the choice of $x_t$ is simply the point of the confidence ellipsoid $B_t$ with largest norm. Note also that we present here a simplified variant where the temporal horizon $n$ is known: the original version of the algorithm is anytime.
We now recall Theorem 2 of \citep{dani2008stochastic}.

\begin{theorem}[$ConfidenceBall_2$]\label{th:dani}

Assume that $(\eta_t)$ is an i.i.d.~white noise, independent of the $(x_{t'})_{t'\leq t}$ and that for all $k = \{1,\ldots,d\}$, $\exists \sigma_k$ such that for all $t$, $|\eta_{t,k}| \leq \frac{1}{2}\sigma_k$.
For large enough $n$, we have with probability $1-\delta$ the following bound for the regret of $ConfidenceBall_2(\D_{d},\delta)$:

\begin{equation*}
R_n(\Alg_{CB_2}) \leq 64 d\Big(||\theta||_2 + ||\si||_2\Big)(\log(n^2/\delta))^2\sqrt{n}.
\end{equation*}

%
%
 
\end{theorem}

\section{The algorithm SL-UCB}\label{s:SL-UCB}

Now we come back to our setting where $n \ll K$.
We present here an algorithm, called {\em Sparse Linear Upper Confidence Bound} (SL-UCB).

\subsection{Presentation of the algorithm}

SL-UCB is divided in two main parts, (i) a first non-adaptive phase, that uses an idea from Compressed Sensing, which is referred to as \textit{support exploration phase} where we project $\theta$ on isotropic random vectors in order to select the arms that belong to what we call the \textit{active set} $\A$, and (ii) a second phase that we call \textit{restricted linear bandit phase} where we apply a linear bandit algorithm to the active set $\A$ in order to balance exploration and exploitation and further minimize the regret. Note that the length of the support exploration phase is problem dependent.

This algorithm takes as parameters: $\bar{\sigma}_2$ and $\bar{\theta}_2$ which are upper bounds respectively on $||\sigma||_2$ and $||\theta||_2$, and $\delta$ which is a (small) probability.

First, we define an \textit{exploring set} as 
\vspace{-0.2cm}
\begin{equation}\label{eq:exploring}
 \Er = \frac{1}{\sqrt{K}}\{-1,+1\}^K.
\end{equation}

Note that $\Er \subset \D_{K}$. We sample this set uniformly during the support exploration phase. This gives us some insight about the directions on which the parameter $\theta$ is sparse, using very simple concentration tools\footnote{Note that this idea is very similar to the one of Compressed Sensing.}: at the end of this phase, the algorithm selects a set of coordinates $\A$, named \textit{active set}, which are the directions where $\theta$ is likely to be non-zero. 
The algorithm automatically adapts the length of this phase and that no knowledge of $||\theta||_{2}$ is required. The Support Exploration Phase ends at the first time $t$ such that (i) $\max_k |\hat{\theta}_{k,t}| - \frac{2b}{\sqrt{t}} \geq 0$ for a well-defined constant $b$ and (ii) $t\geq \frac{\sqrt{n}}{\max_k |\hat{\theta}_{k,t}| - \frac{b}{\sqrt{t}}}$.

We then exploit the information collected in the first phase, i.e.~the active set $\A$, by playing a linear bandit algorithm on the intersection of the unit ball $B_{K}$ and the vector subspace spanned by the active set $\A$, i.e.~$Vec(\A)$. Here we choose to use the algorithm $CB_2$ described in \citep{dani2008stochastic}. See Subsection \ref{ss:dani} for an adaptation of this algorithm to our specific case: the set of arms is indeed the unit ball for the $||.||_2$ norm in the vector subspace $Vec(\A)$.

The algorithm is described in Figure~\ref{fig:SL-UCB}.

\begin{figure}[h]
\bookbox{
\begin{algorithmic}
\STATE \textbf{Input:} parameters $\bar{\sigma}_2$, $\bar{\theta}_2$,$\delta$.
\STATE \textbf{Initialize:} Set $b = (\bar{\theta}_2+ \bar{\sigma}_2) \sqrt{2\log(2K/\delta)}$.
\STATE Pull randomly an arm $x_1$ in $\Er$ (defined in Equation~\ref{eq:exploring}) and observe $r_1$
\STATE \textbf{Support Exploration Phase:}
\WHILE{(i) $\max_k |\hat{\theta}_{k,t}| - \frac{2b}{\sqrt{t}} <0$ or (ii) $t < \frac{\sqrt{n}}{\max_k |\hat{\theta}_{k,t}| - \frac{b}{\sqrt{t}}}$}
\STATE Pull randomly an arm $x_t$ in $\Er$ (defined in Equation~\ref{eq:exploring}) and observe $r_t$
\STATE Compute $\hat{\theta}_{t}$ using Equation~\ref{eq:theta}
\STATE Set $t\leftarrow t+1$
\ENDWHILE
\STATE Call $T$ the length of the Support Exploration Phase
\STATE Set $\A = \Big\{k: \hat{\theta}_{k, T} \geq \frac{2b}{\sqrt{T}}\Big\}$
\STATE \textbf{Restricted Linear Bandit Phase:}
  \STATE For $t=T+1,\ldots,n$, apply $CB_2(\D_{K} \cap Vec(\A),\delta)$ and collect the rewards $r_t$.
\end{algorithmic}
}
\caption{The pseudo-code of the SL-UCB algorithm.\label{fig:SL-UCB}}
\end{figure}

Note that the algorithm computes $\hat{\theta}_{k, t}$ using
\vspace{-0.2cm}
\begin{equation}\label{eq:theta}
 \hat{\theta}_{k,t} = \frac{K}{t} \Big( \sum_{i=1}^t x_{k,i} r_i\Big) = \big(\frac{K}{t} X_tR_t\big)_k.
\end{equation}

\vspace{-0.4cm}


\subsection{Main Result}

We first state an assumption on the noise.
\vspace{-0.3cm}
\begin{assumption}\label{as:noise1}
$(\eta_{k,t})_{k,t}$ is an i.i.d.~white noise and $\exists \si_k$ s.t. $|\eta_{k,t}| \leq \frac{1}{2}\si_k$. 
\end{assumption}

Note that this assumption is made for simplicity and that it could easily be generalized to, for instance, sub-Gaussian noise.
Under this assumption, we have the following bound on the regret.

\begin{theorem}\label{thm:m1-regret}
Under Assumption~\ref{as:noise1}, if we choose $\bar{\sigma}_2 \geq ||\si||_2$, and $\bar{\theta}_2 \geq ||\theta||_2$, the regret of SL-UCB is bounded with probability at least $1-5\delta$, as
\vspace{-0.1cm}
\begin{align*}
 &R_n(\Alg_{SL-UCB})\leq 118 (\bar{\theta}_2+ \bar{\sigma}_2)^2 \log(2K/\delta)S \sqrt{n}.
\end{align*}

\end{theorem}

The proof of this result is reported in Section \ref{s:proof}.


The algorithm SL-UCB  first uses an idea of Compressed Sensing: it explores by performing random projections and builds an estimate of $\theta$. It then selects the support as soon as the uncertainty is small enough, and applies $CB_2$ to the selected support. The particularity of this algorithm is that the length of the support exploration phase adjusts to the difficulty of finding the support: the length of this phase is of order $O(\frac{\sqrt{n}}{||\theta||_2})$. More precisely, the smaller $||\theta||_{2}$, the more difficult the problem (since it is difficult to find the largest components of the support), and the longer the support exploration phase. But note that the regret does not deteriorate for small values of $||\theta||_2$ since in such case the loss at each step is small too.

An interesting feature of SL-UCB is that it does not require the knowledge of the sparsity $S$ of the parameter.

\section{The gradient ascent as a bandit problem}\label{s:grad.asc}

The aim of this section is to propose a gradient optimization technique to maximize a function $f: \R^K \rightarrow \R$ \textit{when the dimension $K$ is large compared to the number of gradient steps $n$, i.e.~$n \ll K$}. We assume that the function $f$ depends on a small number of relevant variables: it corresponds to the assumption that the gradient of $f$ is sparse.


We consider a stochastic gradient ascent (see for instance the book of \cite{bertsekas1999nonlinear} for an exhaustive survey on gradient methods), where one estimates the gradient of $f$ at a sequence of points and moves in the direction of the gradient estimate during $n$ iterations. 

%

\subsection{Formalization}


The objective is to apply gradient ascent to a differentiable function $f$ assuming that we are allowed to query this function $n$ times only. We write $u_t$ the $t-$th point where we sample $f$, and choose it such that $||u_{t+1} - u_t||_2 = \epsilon$, where $\epsilon$ is the gradient step.

Note that by the Theorem of intermediate values
\vspace{-0.1cm}
\begin{align*}
f(u_n) - f(u_0) &= \sum_{t=1}^n f(u_t) - f(u_{t-1})\\ 
&= \sum_{t=1}^n \langle (u_{t}-u_{t-1}), \nabla f(w_t)\rangle , 
\end{align*}

where $w_t$ is an appropriate barycenter of $u_t$ and $u_{t-1}$.

\begin{figure*}[!hbtp]
\begin{minipage}{0.45\textwidth}
\includegraphics[width=10cm]{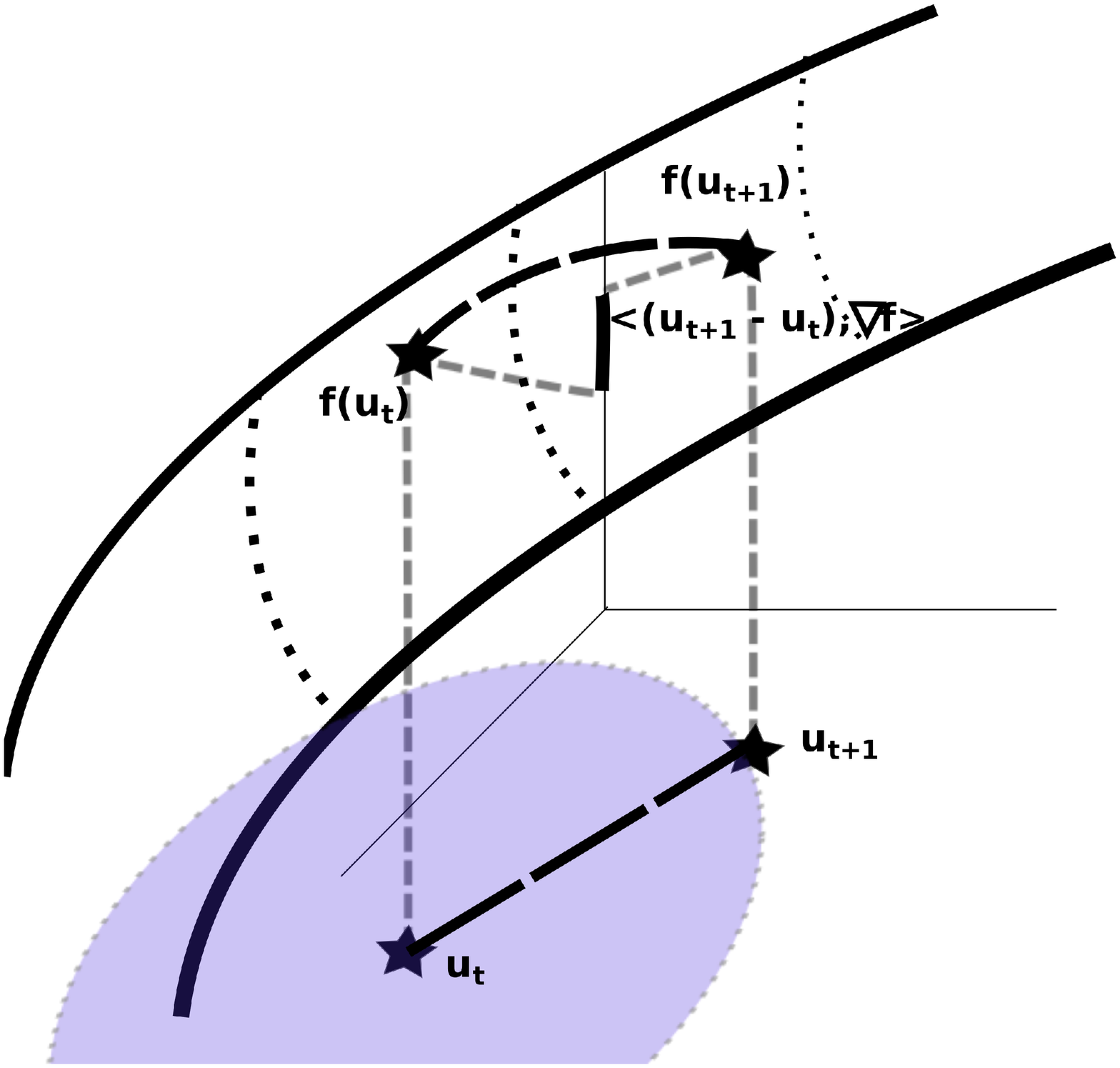}
\end{minipage}
\begin{minipage}{0.45\textwidth}
\includegraphics[width=10cm]{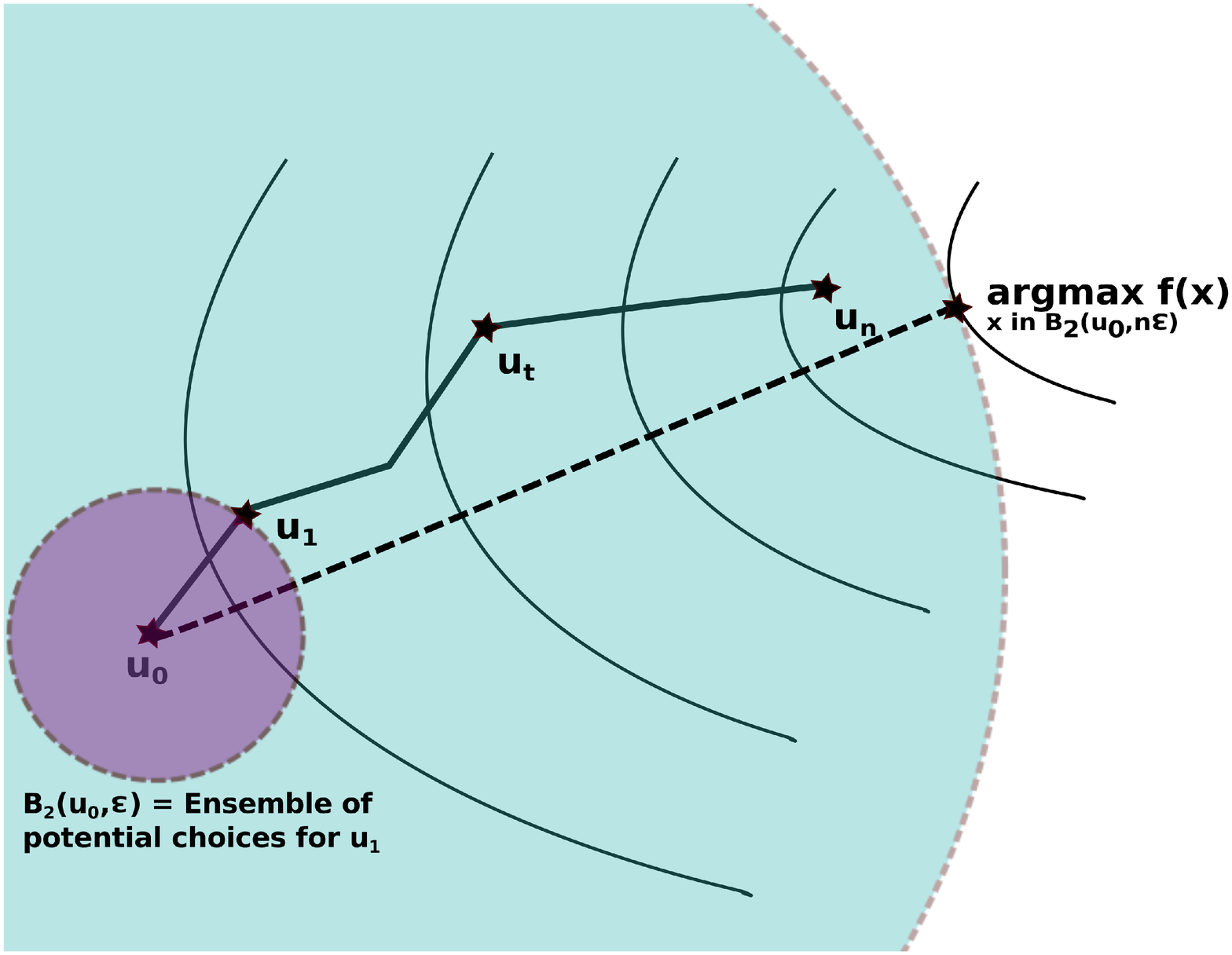}
\end{minipage}
\caption{The gradient ascent: the left picture illustrates the problem written as a linear bandit problem with rewards and the right picture illustrates the regret.}\label{fig:alglinban}
\end{figure*}

We can thus model the problem of  gradient ascent by a linear bandit problem where the reward is what we gain/loose by moving from point $u_{t-1}$ to point $u_{t}$, i.e.~$f(u_t) - f(u_{t-1})$. More precisely, rewriting this problem with previous notations, we  have $\theta + \eta_t = \nabla f(w_t)$\footnote{Note that in order for the model in Section \ref{s:setting} to hold, we need to relax the assumption that $\eta$ is i.i.d..}, and $x_t = u_t - u_{t-1}$. We illustrate this model in Figure~\ref{fig:alglinban}.

If we assume that the function $f$ is (locally) linear and that there are some i.i.d.~measurement errors, we are exactly in the setting of Section~\ref{s:setting}. The objective of minimizing the regret, i.e.,
\begin{equation*}
R_n(\Alg) = \max_{x\in \B_2(u_0,n\epsilon)} f(x) - f(u_n),
\end{equation*}
thus corresponds to the problem of maximizing $f(u_n)$, the $n$-th evaluation of $f$. Thus the regret corresponds to the evaluation of $f$ at the $n$-th step compared to an ideal gradient ascent (that assumes that the true gradient is known and followed for $n$ steps). Applying SL-UCB algorithm implies that the regret is in $O(S\epsilon\sqrt{n})$.

\paragraph{Remark on the noise:}
Assumption~\ref{as:noise1}, which states that the noise added to the function is of the form $\langle u_t - u_{t-1}, \eta_t\rangle $ is specially suitable for gradient ascent because it corresponds to the cases where the noise is an approximation error and depends on the gradient step.

%

\paragraph{Remark on the linearity assumption:}

Matching the stochastic bandit model in Section \ref{s:setting} to the problem of gradient ascent corresponds to assuming that the function is (locally) linear in a neighborhood of $u_0$, and that we have in this neighborhood $f(u_{t+1}) - f(u_t)  = \langle u_{t+1} - u_t, \nabla f(u_0) + \eta_{t+1}\rangle$, where the noise $\eta_{t+1}$ is i.i.d. This setting is somehow restrictive: we made it in order to offer a first, simple solution for the problem. When the function is not linear, one should also consider the additional approximation error.

\subsection{Numerical experiment}

\begin{figure}[!hbtp]
\begin{center}
\includegraphics[width=8.6cm,height=8cm]{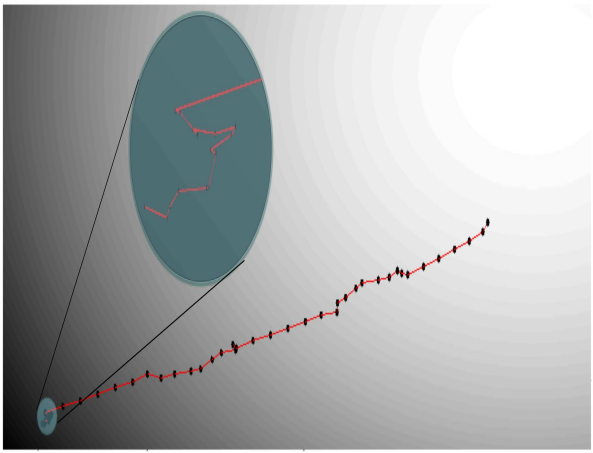}
\end{center}
\caption{Illustration of the trajectory of algorithm SL-UCB with a budget $n=50$, with a zoom at the beginning of the trajectory to illustrate the support exploration phase. The levels of gray correspond to the contours of the function.}
\label{fig:alglinban2}
\end{figure}

In order to illustrate the mechanism of our algorithm, we apply SL-UCB to a quadratic function in dimension $100$ where only two dimensions are informative. Figure~\ref{fig:alglinban2} shows with grey levels the projection of the function onto these two informative directions and a trajectory followed by $n=50$ steps of gradient ascent. The beginning of the trajectory shows an erratic behavior (see the zoom) due to the initial support exploration phase (the projection of the gradient steps onto the  relevant directions are small and random). However, the algorithm quickly selects the righ support of the gradient and the restricted linear bandit phase enables to follow very efficiently the gradient along the two relevant directions.

We now want to illustrate the performances of SL-UCB on more complex problems. We fix the number of pulls to $n=100$, and we try different values of $K$, in order to produce results for different values of the ratio $\frac{K}{n}$. The larger this ratio, the more difficult the problem. We choose a quadratic function that is not constant in $S=10$ directions\footnote{We keep the same function for different values of $K$. It is the quadratic function $f(x) = \sum_{k=1}^{10} -20(x_k - 25)^2$.}. 

We compare our algorithm SL-UCB to two strategies: the {\em ``oracle'' gradient strategy} (OGS), i.e.~a gradient algorithm with access to the \textit{full} gradient of the function\footnote{Each of the $100$ pulls corresponds to an access to the full gradient of the function at a chosen point.}, and the {\em random best direction} (BRD) strategy (i.e., at a given point, chooses a random direction, observes the value of the function a step further in this direction, and moves to that point if the value of the function at this point is larger than its value at the previous point). In Figure~\ref{tab:res}, we report the difference between the value at the final point of the algorithm and the value at the beginning.

\begin{figure}[!hbtp]
\begin{center}

\begin{tabular}{|c||l|l|l|}
  \hline
$K/n$  & OGS & SL-UCB & BRD  \\
\hline
\hline
  2 & $1.875\;10^{5}$ & $1.723\;10^{5}$&  $2.934\;10^{4}$  \\
  \hline
  10 & $1.875\;10^{5}$ & $1.657\;10^{5}$&   $1.335\;10^{4}$ \\
  \hline
  100 & $1.875\;10^{5}$ & $1.552\;10^{5}$  & $5.675\;10^{3}$ \\
  \hline
\end{tabular}
\end{center}
\vspace{-0.5cm}
\caption{We report, for different values of $\frac{K}{n}$ and different strategies, the value of $f(u_n) - f(u_0)$.}\label{tab:res}
\end{figure}

The performances of SL-UCB is (slightly) worse than the optimal ``oracle'' gradient strategy. This is due to the fact that SL-UCB is only given a partial information on the gradient. However it performs much better than the random best direction. Note that the larger $\frac{K}{n}$, the more important the improvements of SL-UCB over the random best direction strategy. This can be explained by the fact that the larger $\frac{K}{n}$, the less probable it is that the random direction strategy picks a direction of interest, whereas our algorithm is designed for efficiently selecting the relevant directions.


\section{Analysis of the SL-UCB algorithm}\label{s:proof}

\subsection{Definition of a high-probability event $\xi$}

\paragraph{Step 0: Bound on the variations of $\hat \theta_t$ around its mean during the Support Exploration Phase}

Note that since $x_{k,t} = \frac{1}{\sqrt{K}}$ or $x_{k,t} = -\frac{1}{\sqrt{K}}$ during the Support Exploration Phase, the estimate $\hat \theta_t$ of $\theta$ during this phase is such that, for any $t_0 \leq T$ and any k

\vspace{-0.3cm}
\begin{eqnarray}
 \hat{\theta}_{k,t_0} &=& \frac{K}{t_0} \Big( \sum_{t=1}^{t_0} x_{k,t} r_t\Big) \nonumber\\
 &=&  \frac{K}{t_0} \Big( \sum_{t=1}^{t_0} x_{k,t} \sum_{k'=1}^K x_{k',t} (\theta_{k'} + \eta_{k',t})\Big) \nonumber \\
&=& \frac{K}{t_0} \sum_{t=1}^{t_0} x_{k,t}^2 \theta_{k} + \frac{K}{t_0} \sum_{t=1}^{t_0} x_{k,t} \sum_{k'\neq k} x_{k',t} \theta_{k'} \nonumber\\ 
& &+ \frac{K}{t_0} \sum_{t=1}^{t_0} x_{k,t} \sum_{k'=1}^K x_{k',t} \eta_{k',t} \nonumber \\
&=& \theta_{k} + \frac{1}{t_0} \sum_{t=1}^{t_0} \sum_{k'\neq k} b_{k,k',t} \theta_{k'} \nonumber\\ 
& &+ \frac{1}{t_0} \sum_{t=1}^{t_0}  \sum_{k'=1}^K b_{k,k',t} \eta_{k',t},\label{eq:exprthet}
\end{eqnarray}
\vspace{-0.3cm}

where $b_{k,k',t} = K x_{k,t}x_{k',t}$.

Note that since the $x_{k,t}$ are i.i.d.~random variables such that $x_{k,t} = \frac{1}{\sqrt{K}}$ with probability $1/2$ and $x_{k,t} = -\frac{1}{\sqrt{K}}$ with probability $1/2$, the $(b_{k,k',t})_{k'\neq k,t}$ are i.i.d.~Rademacher random variables, and $b_{k,k,t} = 1$.

\paragraph{Step 1: Study of the first term.} Let us first study $\frac{1}{t_0} \sum_{t=1}^{t_0} \sum_{k'\neq k} b_{k,k',t} \theta_{k'}$.

Note that the $b_{k,k',t} \theta_{k'}$ are $(K-1)T$ zero-mean independent random variables and that among them, $\forall k' \in \{1,...,K\}$, $t_0$ of them are bounded by $\theta_{k'}$, i.e.~the $(b_{k,k',t} \theta_{k'})_t$. By Hoeffding's inequality, we thus have with probability $1-\delta$ that $|\frac{1}{t_0} \sum_{t=1}^{t_0} \sum_{k'\neq k}^K b_{k,k',t} \theta_{k'}| \leq \frac{||\theta||_2 \sqrt{2\log(2/\delta)}}{\sqrt{t_0}}$. Now by using an union bound on all the $k=\{1,\ldots,K\}$, we have w.p.~$1-\delta$, $\forall k$,

\vspace{-0.3cm}
\begin{equation}\label{eq:part1}
|\frac{1}{t_0} \sum_{t=1}^{t_0} \sum_{k'\neq k} b_{k,k',t} \theta_{k'}| \leq \frac{||\theta||_2 \sqrt{2\log(2K/\delta)}}{\sqrt{t_0}}.
\end{equation}

%
%
%
%
%
%

\paragraph{Step 2: Study of the second term.} Let us now study $\frac{1}{t_0} \sum_{t=1}^{t_0}  \sum_{k'=1}^K b_{k,k',t} \eta_{k',t}$.

Note that the $(b_{k,k',t} \eta_{k',t})_{k',t}$ are $Kt_0$ independent zero-mean random variables, and that among these variables, $\forall k \in \{1,...,K\}$, $t_0$ of them are bounded by $\frac{1}{2}\sigma_k$. By Hoeffding's inequality, we thus have with probability $1-\delta$, $|\frac{1}{t_0} \sum_{t=1}^{t_0} \sum_{k'=1}^K b_{k,k',t} \eta_{k',t}| \leq \frac{||\sigma||_2 \sqrt{2\log(2/\delta)}}{\sqrt{t_0}}$. Thus by an union bound, with probability $1-\delta$, $\forall k$,

\vspace{-0.3cm}

\begin{equation}\label{eq:part2}
|\frac{1}{T} \sum_{t=1}^{t_0} \sum_{k'=1}^K b_{k,k',t} \eta_{k',t}| \leq \frac{||\sigma||_2 \sqrt{2\log(2K/\delta)}}{\sqrt{t_0}}.
\end{equation}

\paragraph{Step 3: Final bound.} Finally for a given $t_0$, with probability $1-2\delta$, we have by Equations~\ref{eq:exprthet}, \ref{eq:part1} and \ref{eq:part2}

\vspace{-0.3cm}
\begin{equation}\label{eq:part3}
 ||\hat{\theta}_{T} - \theta||_{\infty} \leq \frac{(||\theta||_2+ ||\sigma||_2) \sqrt{2\log(2K/\delta)}}{\sqrt{T}}.
\end{equation}

\paragraph{Step 4: Definition of the event of interest.}

Now we consider the event $\xi$ such that

\vspace{-0.3cm}
\begin{align}
 \xi &= \bigcap_{t=1,\ldots,n} \Bigg\{ \omega \in \Omega / ||\theta - \frac{K}{t} X_{t}R_{t} ||_{\infty} \leq  \frac{b}{\sqrt{t}} \Bigg\}, \label{eq:def.xi}
\end{align}

where $b=(\bar{\theta}_2+ \bar{\sigma}_2) \sqrt{2\log(2K/\delta)}$.

From Equation \ref{eq:part3} and an union bound over time, we deduce that $\P(\xi) \geq 1-2n\delta$.

\vspace{-0.3cm}

\subsection{Length of the Support Exploration Phase}\label{ss:length}

The Support Exploration Phase ends at the first time $t$ such that (i) $\max_k |\hat{\theta}_{k,t}| - \frac{2b}{\sqrt{t}} >0$ and (ii) $t\geq \frac{\sqrt{n}}{\max_k |\hat{\theta}_{k,t}| - \frac{b}{\sqrt{t}}}$.

\paragraph{Step 1: A result on the empirical best arm}

On the event $\xi$, we know that for any $t$ and any $k$, $|\theta_{k}| - \frac{b}{\sqrt{t}} \leq |\hat{\theta}_{k,t}| \leq |\theta_{k}| + \frac{b}{\sqrt{t}}$. In particular for $k^* = \arg\max_k |\theta_k|$ we have
\begin{align}\label{eq:criterion}
|\theta_{k^*}| - \frac{b}{\sqrt{t}} \leq \max_k |\hat{\theta}_{k,t}| \leq |\theta_{k^*}| + \frac{b}{\sqrt{t}}.
\end{align}

\vspace{-0.3cm}

\paragraph{Step 2: Maximum length of the Support Exploration Phase.}

If $|\theta_{k^*}| - \frac{3b}{\sqrt{t}} >0$ then by Equation~\ref{eq:criterion}, the first (i) criterion is verified on $\xi$. If $t \geq \frac{1}{\theta_{k^*} - \frac{3b}{\sqrt{t}}}\sqrt{n}$ then by Equation~\ref{eq:criterion}, the second (ii) criterion is verified on $\xi$.

Note that both those conditions are thus verified if $t \geq \max\big(\frac{9b^2}{|\theta_{k^*}|^2}, \frac{4\sqrt{n}}{3|\theta_{k^*}|}\big)$. The Support Exploration Phase stops thus before this moment. Note that as the budget of the algorithm is $n$, we have on $\xi$ that $T \leq \max\big(\frac{9b^2}{|\theta_{k^*}|^2}, \frac{4\sqrt{n}}{3|\theta_{k^*}|}, n\big) \leq \frac{9\sqrt{S}b^2}{||\theta||_2} \sqrt{n}$. We write $T_{\max} = \frac{9\sqrt{S}b^2}{||\theta||_2} \sqrt{n}$.

\paragraph{Step 3: Minimum length of the Support Exploration Phase.}

If the first (i) criterion is verified then on $\xi$ by Equation~\ref{eq:criterion} $|\theta_{k^*}| - \frac{b}{\sqrt{t}} >0$. If the second (ii) criterion is verified then on $\xi$ by Equation~\ref{eq:criterion} we have $t \geq \frac{\sqrt{n}}{|\theta_{k^*}|}$.

Combining those two results, we have on the event $\xi$ that $T \geq \max\big(\frac{b^2}{\theta_{k^*}^2}, \frac{\sqrt{n}}{|\theta_{k^*}|}\big) \geq \frac{b^2}{||\theta||_2}\sqrt{n}$. We write $T_{\min} = \frac{b^2}{||\theta||_2} \sqrt{n}$.

\subsection{Description of the set $\A$}

The set $\A$ is defined as $\A = \Big\{k: |\hat{\theta}_{k, T}| \geq \frac{2b}{\sqrt{T}}\Big\}$.

\paragraph{Step 1: Arms that are in $\A$}

Let us consider an arm $k$ such that $|\theta_k| \geq \frac{3b\sqrt{||\theta||_2}}{n^{1/4}}$. Note that $T \geq T_{\min} = \frac{b^2}{||\theta||_2} \sqrt{n}$ on $\xi$. We thus know that on $\xi$
\vspace{-0.3cm}
\begin{align*}
|\hat{\theta}_{k,T}| &\geq |\theta_{k}| - \frac{b}{\sqrt{T}}\geq  \frac{3b\sqrt{||\theta||_2}}{n^{1/4}} -  \frac{b\sqrt{||\theta||_2}}{n^{1/4}} \geq \frac{2b}{\sqrt{T}}.
\end{align*}
This means that $k \in \A$ on $\xi$. We thus know that $|\theta_k| \geq \frac{3b\sqrt{||\theta||_2}}{n^{1/4}}$ implies on $\xi$ that $k \in \A$.

\paragraph{Step 2: Arms that are not in $\A$}

Now let us consider an arm $k$ such that $|\theta_k| < \frac{b}{2\sqrt{n}}$. Then on $\xi$, we know that
\vspace{-0.3cm}
\begin{align*}
|\hat{\theta}_{k,T}| &< |\theta_{k}| + \frac{b}{\sqrt{T}}< \frac{b}{2\sqrt{n}} + \frac{b}{\sqrt{T}} < \frac{3b}{2\sqrt{T}} < \frac{2b}{\sqrt{T}}.
\end{align*}

This means that $k \in \A^c$ on $\xi$. This implies that on $\xi$, if $|\theta_k|=0$, then $k \in \A^c$.

\paragraph{Step 3: Summary.}

Finally, we know that $\A$ is composed of all the $|\theta_k| \geq \frac{3b\sqrt{||\theta||_2}}{n^{1/4}}$, and that it contains only the strictly positive components $\theta_k$, i.e.~at most $S$ elements since $\theta$ is $S-$sparse. We write $\A_{\min} = \{ k: |\theta_k| \geq \frac{3b\sqrt{||\theta||_2}}{n^{1/4}}\}$.


\subsection{Comparison of the best element on $\A$ and on $\D_{K}$.}

Now let us compare $\max_{x_t \in Vec(\A)\cap\D_{K}} \langle \theta,x_t\rangle$ and  $\max_{x_t \in \D_{K}} \langle \theta,x_t\rangle$.

At first, note that $\max_{x_t \in \D_{K}} \langle \theta,x_t\rangle = ||\theta||_2$ and that $\max_{x_t \in Vec(\A)\cap\D_{K}} \langle \theta,x_t\rangle = ||\theta_{\A}||_2 = \sqrt{\sum_{k=1}^K \theta_k^2 \ind{k \in \A}}$, where $\theta_{\A,k} = \theta_k$ if $k \in \A$ and $\theta_{\A,k} = 0$ otherwise. This means that
\vspace{-0.2cm}
\begin{align}
 &\max_{x_t \in \D_{K}} \langle \theta,x_t\rangle - \max_{x_t \in Vec(\A)\cap\D_{K}} \langle \theta,x_t\rangle \nonumber \\ 
&= ||\theta||_2 - ||\theta\ind{k\in\A}||_2 = \frac{||\theta||_2^2 - ||\theta\ind{k\in\A}||_2^2}{||\theta||_2 + ||\theta\ind{k\in\A}||_2} \nonumber \\
&\leq \frac{\sum_{k\in\A^c} \theta_k^2}{||\theta||_2}\leq \frac{\sum_{k\in\A_{\min}^c} \theta_k^2}{||\theta||_2} \leq \frac{9Sb^2}{\sqrt{n}}. \label{eq:distthetbig}
\end{align}
\vspace{-0.4cm}

\subsection{Expression of the regret of the algorithm}

Assume that we run the algorithm $CB_2(Vec(\A)\cap\D_{K},\delta, T)$ at time $T$ where $\A \subset Supp(\theta)$ with a budget of $n_1 = n-T$ samples. In the paper~\citep{dani2008stochastic}, they prove that on an event $\xi_2(Vec(\A)\cap\D_{K}, \delta, T)$ of probability $1-\delta$ the regret of algorithm $CB_2$ is bounded by $R_n(\Alg_{CB_2(Vec(\A)\cap\D_{K},\delta, T)}) \leq 64 |\A|\Big(||\theta||_2 + ||\si||_2\Big)(\log(n^2/\delta))^2\sqrt{n_1}$.

Note that since $\A \subset Supp(\theta)$, we have $\xi_2(Vec(\A)\cap\D_{K}, \delta, T) \subset \xi_2(Vec(Supp(\theta))\cap\D_{K}, \delta, T)$ (see the paper~\citep{dani2008stochastic} for more details on the event $\xi_2$). We thus now that, conditionally to $T$, with probability $1-\delta$, the regret is bounded for any $\A \subset Supp(\theta)$ as $R_n(\Alg_{CB_2(Vec(\A)\cap\D_{K},\delta,T)}) \leq 64 S\Big(||\theta||_2 + ||\si||_2\Big)(\log(n^2/\delta))^2\sqrt{n_1}$.

By an union bound on all possible values for $T$ (i.e.~from $1$ to $n$), we obtain that on an event $\xi_2$ whose probability is larger than $1-\delta$, $R_n(\Alg_{CB_2(Vec(\A)\cap\D_{K},\delta, T)}) \leq 64 S\Big(||\theta||_2 + ||\si||_2\Big)(\log(n^3/\delta))^2\sqrt{n}$.

We thus have on $\xi \bigcup \xi_2$, i.e.~on an event with probability larger than $1-2\delta$, that

\vspace{-0.4cm}
\begin{align}
 R_n(\Alg_{SL-UCB}, \delta) &\leq 2T_{\max}||\theta||_2 \nonumber\\ 
&+ \max_t R_n(\Alg_{CB_2(Vec(\A)\cap\D_{K},\delta, t)}) \nonumber\\
  &+ n \Big( \max_{x \in \D_{K}} \langle x,\theta\rangle  - \hspace{-6mm}\max_{x \in \D_{K}\cap Vect(\A_{\min})} \langle x,\theta\rangle \Big). \nonumber
\end{align}
\vspace{-0.4cm}

By using this Equation, the maximal length of the support exploration phase $T_{\max}$ deduced in Step 2 of Subsection~\ref{ss:length}, and Equation~\ref{eq:distthetbig}, we obtain on $\xi$ that
\begin{eqnarray*}
 R_n &\leq& 64 S\big(||\theta||_2 + ||\si||_2\big)(\log(n^2/\delta))^2\sqrt{n} \\
& & + 18 Sb^2\sqrt{n}+ 9Sb^2 \sqrt{n}\\
&\leq& 118 (\bar{\theta}_2+ \bar{\sigma}_2)^2 \log(2K/\delta)S \sqrt{n}.
\end{eqnarray*}
by using $b=(\bar{\theta}_2+ \bar{\sigma}_2) \sqrt{2\log(2K/\delta)}$ for the third step.

\section*{Conclusion}

In this paper we introduced the SL-UCB algorithm for sparse linear bandits in high dimension. It has been designed using ideas from Compressed Sensing and Bandit Theory. Compressed Sensing is used in the support exploration phase, in order to select the support of the parameter. A linear bandit algorithm is then applied to the small dimensional subspace defined in the first phase. We derived a regret bound of order $O(S\sqrt{n})$. Note that the bound scales with the sparsity $S$ of the unknown parameter $\theta$ instead of the dimension $K$ of the parameter (as is usually the case in linear bandits). We then provided an example of application for this setting, the optimization of a function in high dimension. Possible further research directions include:
\vspace{-0.2cm}
\begin{itemize}
\item The case when the support of $\theta$ changes with time, for which it would be important to define assumptions under which sub-linear regret is achievable. One idea would be to use techniques developed for \textit{adversarial bandits} (see \citep{abernethy2008competing, bartlett2008high, cesa2009combinatorial, koolen2010hedging, audibert2011minimax}, but also \citep{flaxman2005online} for a more gradient-specific modeling) or also from \textit{restless/switching bandits} (see e.g.~\citep{whittle1988restless,nino2001restless, slivkins2008adapting,garivier} and many others). This would be particularly interesting to model gradient ascent for e.g.~convex function where the support of the gradient is not constant.
\item Designing an improved analysis (or algorithm) in order to achieve a regret of order $O(\sqrt{Sn})$, which is the lower bound for the problem of linear bandits in a space of dimension $S$. Note that when an upper bound $S'$ on the sparsity is available, it seems possible to obtain such a regret by replacing condition (ii) in the algorithm by $t < \frac{\sqrt{n}}{||\big(\hat{\theta}_{t,k}\ind{\hat{\theta}_{t,k} \geq \frac{b}{\sqrt{t}}}\big)_k||_2 - \frac{\sqrt{S'b}}{\sqrt{t}}}$, and using for the Exploitation phase the algorithm in~\citep{rusmevichientong2008linearly}. The regret of such an algorithm would be in $O(\sqrt{S'n})$. But it is not clear whether it is possible to obtain such a result when no upper bound on $S$ is available (as is the case for SL-UCB).
\end{itemize}

\subsection*{Acknowledgements}

\vspace{-3mm}
This research was partially supported by Region Nord-Pas-de-Calais Regional Council, French ANR EXPLO-RA (ANR-08-COSI-004), the European Community’s Seventh Framework Programme (FP7/2007-2013) under grant agreement 231495 (project CompLACS), and by Pascal-2.

\newpage
\bibliography{allocation.bib}

\begin{thebibliography}{19}
\providecommand{\natexlab}[1]{#1}
\providecommand{\url}[1]{\texttt{#1}}
\expandafter\ifx\csname urlstyle\endcsname\relax
  \providecommand{\doi}[1]{doi: #1}\else
  \providecommand{\doi}{doi: \begingroup \urlstyle{rm}\Url}\fi

\bibitem[A.~Garivier(2011)]{garivier}
E.~Moulines A.~Garivier.
\newblock On upper-confidence bound policies for non-stationary bandit
  problems.
\newblock In \emph{Algorithmic Learning Theory (ALT)}, 2011.

\bibitem[Abbasi-yadkori et~al.(2011)Abbasi-yadkori, Pal, and Szepesvari]{Csaba}
Y.~Abbasi-yadkori, D.~Pal, and C.~Szepesvari.
\newblock Improved algorithms for linear stochastic bandits.
\newblock In \emph{Advances in Neural Information Processing Systems}, 2011.

\bibitem[Abbasi-yadkori et~al.(2012)Abbasi-yadkori, Pal, and Szepesvari]{Yasin}
Y.~Abbasi-yadkori, D.~Pal, and C.~Szepesvari.
\newblock Online-to-confidence-set conversions and application to sparse
  stochastic bandits.
\newblock In \emph{Artificial Intelligence and Statistics}, 2012.

\bibitem[Abernethy et~al.(2008)Abernethy, Hazan, and
  Rakhlin]{abernethy2008competing}
J.~Abernethy, E.~Hazan, and A.~Rakhlin.
\newblock Competing in the dark: An efficient algorithm for bandit linear
  optimization.
\newblock In \emph{Proceedings of the 21st Annual Conference on Learning Theory
  (COLT)}, volume~3. Citeseer, 2008.

\bibitem[Audibert et~al.(2011)Audibert, Bubeck, and
  Lugosi]{audibert2011minimax}
J.Y. Audibert, S.~Bubeck, and G.~Lugosi.
\newblock Minimax policies for combinatorial prediction games.
\newblock \emph{Arxiv preprint arXiv:1105.4871}, 2011.

\bibitem[Bartlett et~al.(2008)Bartlett, Dani, Hayes, Kakade, Rakhlin, and
  Tewari]{bartlett2008high}
P.L. Bartlett, V.~Dani, T.~Hayes, S.M. Kakade, A.~Rakhlin, and A.~Tewari.
\newblock High-probability regret bounds for bandit online linear optimization.
\newblock In \emph{Proceedings of the 21st Annual Conference on Learning Theory
  (COLT 2008)}, pages 335--342. Citeseer, 2008.

\bibitem[Bertsekas(1999)]{bertsekas1999nonlinear}
D.P. Bertsekas.
\newblock \emph{Nonlinear programming}.
\newblock Athena Scientific Belmont, MA, 1999.

\bibitem[Blumensath and Davies(2009)]{blumensath2009iterative}
T.~Blumensath and M.E. Davies.
\newblock Iterative hard thresholding for compressed sensing.
\newblock \emph{Applied and Computational Harmonic Analysis}, 27\penalty0
  (3):\penalty0 265--274, 2009.

\bibitem[Candes and Tao(2007)]{candes2007dantzig}
E.~Candes and T.~Tao.
\newblock The dantzig selector: statistical estimation when p is much larger
  than n.
\newblock \emph{The Annals of Statistics}, 35\penalty0 (6):\penalty0
  2313--2351, 2007.

\bibitem[Cesa-Bianchi and Lugosi(2009)]{cesa2009combinatorial}
N.~Cesa-Bianchi and G.~Lugosi.
\newblock Combinatorial bandits.
\newblock In \emph{Proceedings of the 22nd Annual Conference on Learning Theory
  (COLT 09)}. Citeseer, 2009.

\bibitem[Chen et~al.(1999)Chen, Donoho, and Saunders]{chen1999atomic}
S.S. Chen, D.L. Donoho, and M.A. Saunders.
\newblock Atomic decomposition by basis pursuit.
\newblock \emph{SIAM journal on scientific computing}, 20\penalty0
  (1):\penalty0 33--61, 1999.

\bibitem[Dani et~al.(2008)Dani, Hayes, and Kakade]{dani2008stochastic}
V.~Dani, T.P. Hayes, and S.M. Kakade.
\newblock Stochastic linear optimization under bandit feedback.
\newblock In \emph{Proceedings of the 21st Annual Conference on Learning Theory
  (COLT)}. Citeseer, 2008.

\bibitem[Filippi et~al.(2010)Filippi, Capp{\'e}, Garivier, and
  Szepesv{\'a}ri]{filippiparametric}
S.~Filippi, O.~Capp{\'e}, A.~Garivier, and C.~Szepesv{\'a}ri.
\newblock Parametric bandits: The generalized linear case.
\newblock In \emph{Advances in Neural Information Processing Systems}, 2010.

\bibitem[Flaxman et~al.(2005)Flaxman, Kalai, and McMahan]{flaxman2005online}
A.D. Flaxman, A.T. Kalai, and H.B. McMahan.
\newblock Online convex optimization in the bandit setting: gradient descent
  without a gradient.
\newblock In \emph{Proceedings of the sixteenth annual ACM-SIAM symposium on
  Discrete algorithms}, pages 385--394. Society for Industrial and Applied
  Mathematics, 2005.

\bibitem[Koolen et~al.(2010)Koolen, Warmuth, and Kivinen]{koolen2010hedging}
W.M. Koolen, M.K. Warmuth, and J.~Kivinen.
\newblock Hedging structured concepts.
\newblock In \emph{Proceedings of the 23rd Annual Conference on Learning Theory
  (COLT 19)}. Omnipress, 2010.

\bibitem[Nino-Mora(2001)]{nino2001restless}
J.~Nino-Mora.
\newblock Restless bandits, partial conservation laws and indexability.
\newblock \emph{Advances in Applied Probability}, 33\penalty0 (1):\penalty0
  76--98, 2001.

\bibitem[Rusmevichientong and Tsitsiklis(2008)]{rusmevichientong2008linearly}
P.~Rusmevichientong and J.N. Tsitsiklis.
\newblock Linearly parameterized bandits.
\newblock \emph{Arxiv preprint arXiv:0812.3465}, 2008.

\bibitem[Slivkins and Upfal(2008)]{slivkins2008adapting}
A.~Slivkins and E.~Upfal.
\newblock Adapting to a changing environment: The brownian restless bandits.
\newblock In \emph{Proc. 21st Annual Conference on Learning Theory}, pages
  343--354. Citeseer, 2008.

\bibitem[Whittle(1988)]{whittle1988restless}
P.~Whittle.
\newblock Restless bandits: Activity allocation in a changing world.
\newblock \emph{Journal of applied probability}, pages 287--298, 1988.

\end{thebibliography}
\bibliographystyle{plainnat}
\newpage

\end{document}